\documentclass[a4paper,reqno,11pt]{amsart}

\usepackage{amsmath, amsfonts, amssymb, amsthm, amscd}
\usepackage{graphicx}
\usepackage{psfrag}
\usepackage{perpage}
\usepackage{url}
\usepackage{color}
\usepackage{mathrsfs}
\usepackage{tikz}
\usepackage{mathabx}

\usepackage{dsfont} 

\usepackage[utf8]{inputenc}
\usepackage[T1]{fontenc}
\usepackage{microtype}

\usepackage[a4paper,scale={0.72,0.74},marginratio={1:1},footskip=7mm,headsep=10mm]{geometry}

\usepackage{hyperref}

\makeatletter
\def\@secnumfont{\bfseries\scshape}

\def\section{\@startsection{section}{1}%
  \z@{.7\linespacing\@plus\linespacing}{.5\linespacing}%
  {\normalfont\large\bfseries\scshape\centering}}

\def\subsection{\@startsection{subsection}{2}%
  \z@{.5\linespacing\@plus.7\linespacing}{-.5em}%
  {\normalfont\bfseries\scshape}}

\def\subsubsection{\@startsection{subsubsection}{3}%
  \z@{.5\linespacing\@plus.7\linespacing}{-.5em}%
  {\normalfont\scshape}}

\def\specialsection{\@startsection{section}{1}%
  \z@{\linespacing\@plus\linespacing}{.5\linespacing}%
  {\normalfont\centering\large\bfseries\scshape}}
\makeatother

\makeatletter

\renewenvironment{proof}[1][\proofname]{\par
\pushQED{\qed}%
\normalfont \topsep4\p@\@plus4\p@\relax
\trivlist
\item[\hskip\labelsep
\bfseries
#1\@addpunct{.}]\ignorespaces
}{%
\popQED\endtrivlist\@endpefalse
}
\makeatother

\makeatletter
\newcommand \Dotfill {\leavevmode \leaders \hb@xt@ 6pt{\hss .\hss }\hfill \kern \z@}
\makeatother

\makeatletter
\def\@tocline#1#2#3#4#5#6#7{\relax
  \ifnum #1>\c@tocdepth 
  \else
    \par \addpenalty\@secpenalty\addvspace{#2}%
    \begingroup \hyphenpenalty\@M
    \@ifempty{#4}{%
      \@tempdima\csname r@tocindent\number#1\endcsname\relax
    }{%
      \@tempdima#4\relax
    }%
    \parindent\z@ \leftskip#3\relax \advance\leftskip\@tempdima\relax
    \rightskip\@pnumwidth plus4em \parfillskip-\@pnumwidth
    #5\leavevmode\hskip-\@tempdima
      \ifcase #1
       \or\or \hskip 1.65em \or \hskip 3.3em \else \hskip 4.95em \fi%
      #6\nobreak\relax
    \Dotfill
    \hbox to\@pnumwidth{\@tocpagenum{#7}}\par
    \nobreak
    \endgroup
  \fi}
\makeatother

\makeatletter
\def\l@section{\@tocline{1}{0pt}{1pc}{}{\scshape}}
\renewcommand{\tocsection}[3]{%
\indentlabel{\@ifnotempty{#2}{\ignorespaces#1 #2.\hskip 0.7em}}#3}
\def\l@subsection{\@tocline{2}{0pt}{1pc}{5pc}{}}

\def\l@subsubsection{\@tocline{3}{0pt}{1pc}{7pc}{}}

\makeatother

\numberwithin{equation}{section}

\newtheoremstyle{mytheorem}{.7\linespacing\@plus.3\linespacing}{.7\linespacing\@plus.3\linespacing}%
     {\itshape}
     {}
     {\bfseries}
     {. }
     {0.3ex}
     {\thmname{{\bfseries #1}}\thmnumber{ {\bfseries #2}}\thmnote{ (#3)}}  

\theoremstyle{mytheorem}

\newtheorem{theorem}{Theorem}[section]

\newtheorem{remark}[theorem]{Remark}

\newtheorem{assumption}[theorem]{Assumption}

\newcommand{\bbE}{{\ensuremath{\mathbb E}} }

\newcommand{\bbP}{{\ensuremath{\mathbb P}} }

\newcommand{\bbR}{{\ensuremath{\mathbb R}} }

\newcommand{\bbZ}{{\ensuremath{\mathbb Z}} }

\newcommand{\ga}{\alpha}
\newcommand{\gb}{\beta}

\newcommand{\gd}{\delta}

\newcommand{\gs}{\sigma}

\newcommand{\go}{\omega}

\renewcommand{\tilde}{\widetilde}          
\DeclareMathSymbol{\leqslant}{\mathalpha}{AMSa}{"36} 
\DeclareMathSymbol{\geqslant}{\mathalpha}{AMSa}{"3E} 
\DeclareMathSymbol{\eset}{\mathalpha}{AMSb}{"3F}     

\newcommand{\sumtwo}[2]{\sum_{\substack{#1 \\ #2}}} 

\newcommand{\R}{\mathbb{R}}

\newcommand{\Z}{\mathbb{Z}}
\newcommand{\N}{\mathbb{N}}

\def\bs{\boldsymbol}

\newcommand{\PEfont}{\mathrm}

\newcommand{\p}{\ensuremath{\PEfont P}}

\DeclareMathOperator{\e}{\ensuremath{\PEfont E}}

\newcommand\Ro{\ensuremath{\bs{\mathrm{R}}}}

\newcommand{\E}{\e}
\renewcommand{\P}{\p}

\DeclareMathOperator{\bbvar}{\ensuremath{\mathbb{V}ar}}

\newcommand{\ind}{\mathds{1}}

\renewcommand{\epsilon}{\varepsilon}
\renewcommand{\theta}{\vartheta}
\renewcommand{\rho}{\varrho}

\newenvironment{myenumerate}{%
\renewcommand{\theenumi}{\arabic{enumi}}%
\renewcommand{\labelenumi}{{\rm(\theenumi)}}%
\begin{list}{\labelenumi}
	{%
	\setlength{\itemsep}{0.4em}%
	\setlength{\topsep}{0.5em}%
	\setlength\leftmargin{2.45em}%
	\setlength\labelwidth{2.05em}%
	\setlength{\labelsep}{0.4em}%
	\usecounter{enumi}%
	}%
	}%
{\end{list}
}

{\end{list}
}

{\end{list}
}

{\end{myenumerate}}

\newenvironment{myitemize}{%
\begin{list}{$\bullet$}%
 	{%
	\setlength{\itemsep}{0.4em}%
	\setlength{\topsep}{0.5em}%
	\setlength\leftmargin{2.65em}%
	\setlength\labelwidth{2.65em}%
	\setlength{\labelsep}{0.4em}%
	}%
	}%
{\end{list}}

{\end{myitemize}}

\MakePerPage[2]{footnote} 

\date{\today}

\newcommand\dd{\mathrm{d}}

\newcommand\hbeta{{\hat{\beta}}}
\newcommand\hh{{\hat{h}}}

\newcommand\bi{\boldsymbol{i}}

\newcommand\re{\mathrm{ref}}
\newcommand\eff{\mathrm{eff}}

\newcommand\bpsi{\boldsymbol\psi}

\newcommand\bZ{\boldsymbol Z}

\begin{document}

\title[Scaling limits and disorder relevance]
{Scaling limits of disordered systems and disorder relevance}

\begin{abstract}
We review recent works where we have shown that disorder relevance is closely related to the existence of non-trivial, random continuum limits.
 \end{abstract}

\author[F.Caravenna]{Francesco Caravenna}
\address{Dipartimento di Matematica e Applicazioni\\
 Universit\`a degli Studi di Milano-Bicocca\\
 via Cozzi 55, 20125 Milano, Italy}
\email{francesco.caravenna@unimib.it}

\author[R.Sun]{Rongfeng Sun}
\address{Department of Mathematics\\
National University of Singapore\\
10 Lower Kent Ridge Road, 119076 Singapore
}
\email{matsr@nus.edu.sg}

\author[N.Zygouras]{Nikos Zygouras}
\address{Department of Statistics\\
University of Warwick\\
Coventry CV4 7AL, UK}
\email{N.Zygouras@warwick.ac.uk}


\maketitle

\section{Introduction}
\label{sec:intro}

A basic question in disordered systems is whether arbitrarily small amount of disorder changes the critical 
properties (e.g. shift in the critical curve, critical exponents, etc.) of a statistical mechanics model. 
The Harris criterion \cite{H74} suggests that the answer to this question depends on a suitably defined correlation length exponent $\nu$
of the (pure) statistical mechanics model and the dimension $d$. In particular, any amount of disorder, however small,
is sufficient to change the qualitative properties of the system if $\nu<2/d$, while the properties stay the same for small disorder, when
$\nu>2/d$. The case $\nu=2/d$ is marginal and the Harris criterion is inconclusive. There has been some 
effort to mathematically verify Harris' criterion on a case by case basis, showing that certain critical points are shifted 
\cite{AZ09,DGLT09, G07}, or that there is smoothing of phase transitions \cite{AW90, GT06}. 
Earlier efforts to build a mathematical framework around the Harris criterion considered correlation lengths defined via finite size scaling
\cite{CCFS86, CCFS89}.
We propose a new point of view, which focuses on the existence of 
a non-trivial, random continuum limit when disorder scales to zero in a particular way as a function of the lattice spacing.
Our approach covers also the marginal cases (for the moment of polymer type models)
 and further reveals a phase transition that takes place therein.

Let us set our framework.
We define a pure statistical mechanics model as a law $\P^\re_{\Omega_\gd}$ describing the distribution of a
field of (correlated) variables $\sigma = (\sigma_x)_{x\in\Omega_\delta}$,
 on a sub-lattice of  $\Omega\subset \bbR^d$: 
 \begin{equation*}
	\Omega_\delta := \big(\delta^{a_1} \Z \times \delta^{a_2} \Z \times 
	\cdots \times \delta^{a_d} \Z \big) 
	\cap \Omega , \quad
	\text{with} \,\,\,\, \gd>0,\,\,a_1 = 1 \text{ and } a_2, \ldots, a_d > 0,
\end{equation*}
where we allow for different scaling of the dimensions, in which case the \emph{effective dimension} of the model is $d_\eff := a_1 + \ldots + a_d$.
 We assume that the field $\sigma$ is a binary field, e.g. an occupation $\sigma_x\in\{0,1\}$ or spin $\sigma_x\in\{\pm1\}$ field.
Disorder is modeled 
by a family of  i.i.d.\ random variables $\omega := (\omega_x)_{x\in\Omega_\delta}$  
with zero mean, unit variance, and locally finite exponential moments.
Probability and expectation with respect to $\omega$ will be denoted respectively by $\bbP$ and $\bbE$.
Given $\gb>0$, $h\in\R$ and a $\bbP$-typical realization
of the disorder $\omega$, we define the \emph{disordered model}
as the following probability measure $\P^\omega_{{\Omega_\delta}; \gb,h}$ for the field
$\sigma = (\sigma_x)_{x\in\Omega_\delta}$:
\begin{equation} \label{eq:model}
	\P^\omega_{{\Omega_\delta}; \gb,h} (\dd\sigma)
	:= \frac{e^{\sum_{x\in\Omega_\delta} (\gb \omega_x +h) \sigma_x}}
	{Z^\omega_{{\Omega_\delta}; \gb,h}} \P_{\Omega_\delta}^\re (\dd\sigma)\,,
\end{equation}
where the normalizing constant is the partition function defined by
\begin{equation} \label{eq:partfun}
	Z^\omega_{{\Omega_\delta}; \gb,h} := \E_{\Omega_\delta}^\re
	\big[ e^{\sum_{x\in\Omega_\delta} (\gb \omega_x+h) \sigma_x} \big] \,.
\end{equation}
For every $k\in\N$
we define the \emph{$k$-point correlation function} $\psi_{\Omega_\delta}^{(k)}(x_1, \ldots, x_k)$,
for $x_1, \ldots, x_k \in \Omega$, as follows:
denoting by $x_\delta$ the point in $\Omega_\delta$ closest to $x \in \Omega$,
we set
\begin{equation} \label{eq:corrfun}
	\psi_{\Omega_\delta}^{(k)}(x_1, \ldots, x_k) :=
	\begin{cases}
	\E_{\Omega_\delta}^\re \big[\sigma_{(x_1)_\delta} \, \sigma_{(x_2)_\delta}
	\cdots \sigma_{(x_k)_\delta} \big]
	&, \text{if } (x_i)_\delta \ne (x_j)_\delta \text{ for all } i\ne j  \\
	0 &, \text{otherwise}
	\end{cases} .
\end{equation}
We have defined the $k$-point function on the whole $\Omega$, rather
than only on $\Omega_\delta$, for convenience.
\begin{assumption}\label{ass:main}
For every $k\in\N$, there exist a symmetric function $\bpsi_\Omega^{(k)}: \Omega^k \to \R$
and an exponent $\gamma \in [0,\infty)$ such that
\begin{equation} \label{eq:polyconv}
	(\delta^{-\gamma})^k \, \psi_{\Omega_\delta}^{(k)} (x_1, \ldots, x_k)
	\,\xrightarrow[\,\delta \downarrow 0\,]{}\,
	\bpsi_{\Omega}^{(k)} (x_1, \ldots, x_k) \qquad
       \text{in } \,\,L^2(\Omega^k) \,,
\end{equation}
and furthermore, for some $\epsilon > 0$,
\begin{equation} \label{eq:L2sum}
	\limsup_{\ell\to\infty} \limsup_{\delta\downarrow0}\sum_{k>\ell}
	\frac{(1+\epsilon)^k}{k!} \, \big\| \psi_{\Omega_\delta}^{(k)} \big\|_{L^2(\Omega^k)}^2 = 0 \,.
\end{equation}
\end{assumption}
This assumption suggests that the correlation functions of the pure system decay 
 as a power law (if $\gamma>0$) and moreover that their continuum limit exists. This is
  a natural assumption that one expects to hold for systems with a continuous phase transition at {\it critical temperature} 
and which can be verified in many situations, see the examples later on. 
Let us note that we have implicitly assumed here that $\E_{\Omega_\delta}^\re \big[\sigma_{x_\delta} ]$ converges to zero as 
$\gd\downarrow0$. In situations that this limit is not zero, we need in \eqref{eq:corrfun} to consider the centered field 
$\gs_{x_\gd}-\E_{\Omega_\delta}^\re \big[\sigma_{x_\delta} ]$, instead. 
The case $\gamma=0$ is associated to a first order phase
transition at the critical temperature.
The factor $(1+\epsilon)^k$ in  \eqref{eq:L2sum} is needed when one considers an
 $h$ for which the variables $e^{\gb\go_x+h}-1$ are not mean zero. When these are mean zero, then $\epsilon$ 
  can be taken to be equal to zero in \eqref{eq:L2sum}.
  
Under Assumption~\ref{ass:main}, we can determine the scaling limit of the partition function
in the continuum ($\delta \downarrow 0$) and weak disorder ($\gb,h \to 0$) limit:
Fix $\hbeta \ge 0$, $\hh \in \R$ and define
the \emph{continuum partition function} by
\begin{equation} \label{eq:Zcont}
	\bZ^W_{\Omega; \hbeta, \hh} :=
	\sum_{k=0}^{\infty} \;
	\frac{1}{k!} \,
	\idotsint_{\Omega^k}
	\bpsi^{(k)}_{\Omega}(x_1, \ldots, x_k)
	\, \prod_{i=1}^k \big(\hbeta \, W(\dd x_i) +  \hh \, \dd x_i\big) \,,
\end{equation}
where the $k=0$ term of the sum equals $1$ by definition, and
$W(\cdot)$ denotes white noise on (the bounded open set) 
$\Omega \subseteq \R^d$. The convergence of the series in $L^2$ 
is guaranteed by \eqref{eq:L2sum}.
If \eqref{eq:L2sum} holds only with $\epsilon = 0$, 
the series \eqref{eq:Zcont} converges
in $L^p$ for all $p\in(0,2)$. We have 

\begin{theorem}[\cite{CSZ13}]\label{th:mainhere}
Let Assumption~\ref{ass:main} hold.
For fixed $\hbeta > 0$, $\hh \in \R$, rescale $\gb, h$ as follows:
\begin{equation} \label{eq:scalingpar}
	\gb = \gb_\delta = \hbeta \, \delta^{\,d_\eff/2-\gamma} \,, \qquad \ \
	h = h_\delta = \begin{cases}
	\hh \, \delta^{\,d_\eff-\gamma} - \frac{1}{2} \gb^2 &
	\text{in case $\sigma_x \in \{0,1\}$}\\
	\rule{0pt}{1.3em}\hh \, \delta^{\,d_\eff-\gamma} &
	\text{in case $\sigma_x \in \{-1,1\}$}\\
	\end{cases} .
\end{equation}
Then we have the convergence
in distribution to $\bZ^W_{\Omega; \hbeta, \hh}$ 
as $\delta\downarrow 0$:
\begin{equation*}
	\left.\begin{cases}
	Z_{\Omega_\delta; \gb_\delta,h_\delta}^\omega &
	\text{in case $\sigma_x \in \{0,1\}$}\\
	\rule{0pt}{1.5em}e^{-\frac{1}{2} \hbeta^2\,
	\delta^{-2\gamma}} Z_{\Omega_\delta; \gb_\delta,h_\delta}^\omega &
	\text{in case $\sigma_x \in \{-1,1\}$}\\
	\end{cases}\right\}
	\ \xrightarrow[\,\delta \downarrow 0\,]{(d)} \
	\bZ^W_{\Omega; \hbeta, \hh}.
\end{equation*}
\end{theorem}
where the superscript $(d)$ above the arrow denotes convergence in distribution.

The proof has been given in \cite{CSZ13}, cf. Theorems~2.3 and~2.5. Below we give a sketch of it highlighting the main points.
\begin{proof}[Sketch of the proof of Theorem \ref{th:mainhere}.]
Let us consider the so-called {\it high temperature
expansion} ($|\gb|, |h| \ll 1$) of
the partition function $Z^\omega_{\Omega_\delta; \gb, h}$ and for simplicity let us assume that
 the field takes the values $\sigma_x \in \{0,1\}$.
In this case we can write
\begin{equation}\label{eq:Zstart}
	Z^\omega_{\Omega_\delta; \gb, h} =
	\E_{\Omega_\delta}^\re \Bigg[ \prod_{x \in \Omega_\delta} \big(1 + \eta_x \sigma_x
	\big) \Bigg] \,, \qquad \text{where} \qquad
	\eta_x := e^{(\gb \omega_x + h)} - 1 \,.
\end{equation}
 and expanding the product we have
 \begin{equation} \label{eq:Zstart2}
	Z^\omega_{\Omega_\delta; \gb, h} =
	1 + \sum_{k=1}^{|\Omega_\delta|} \;
	\frac{1}{k!} \,
	\sum_{(x_1, x_2, \ldots, x_k) \in (\Omega_\delta)^k}
	\psi^{(k)}_{\Omega_\delta}(x_1, \ldots, x_k)
	\, \prod_{i=1}^k \eta_{x_i} \,,
\end{equation}
where the $k!$ takes into consideration that we sum over
ordered $k$-tuples $(x_1, \ldots, x_k)$.
In this way the partition function is written as a
multi-linear polynomial of the independent random variables $(\eta_x)_{x\in\Omega_\delta}$
with coefficients
given by the $k$-point correlation function of the reference field.
A Taylor expansion shows that when $\gb, h$ are small we have
\begin{equation} \label{eq:Taylor}
	\bbE[\eta_x] \simeq h + \tfrac{1}{2}\gb^2 =: h' \,, \qquad
	\bbvar[\eta_x] \simeq \gb^2 \,.
\end{equation}

When  $|\gb|, |h| \ll 1$, the distribution of
a polynomial chaos expansion, like the right hand side of \eqref{eq:Zstart2}, 
is asymptotically, in the limit $\gd\downarrow 0$, insensitive with respect to
the distribution of the random variables $(\eta_x)_{x\in\Omega_\delta}$,
as long as mean and variance are kept fixed.
The precise formulation of this loosely stated invariance principle
is given by Theorems~2.6 and~2.8 in \cite{CSZ13},
in the form of a Lindeberg principle and is the key point of the method:
denoting by $(\tilde\omega_x)_{x\in\Omega_\delta}$ a family of i.i.d.
standard Gaussians satisfying \eqref{eq:Taylor}, we can approximate, in the limit $\gd\downarrow 0$,
\begin{equation} \label{eq:Zstart3}
	Z^\omega_{\Omega_\delta; \gb, h} \simeq
	1 + \sum_{k=1}^{|\Omega_\delta|} \;
	\frac{1}{k!} \,
	\sum_{(x_1, x_2, \ldots, x_k) \in (\Omega_\delta)^k}
	\psi^{(k)}_{\Omega_\delta}(x_1, \ldots, x_k)
	\, \prod_{i=1}^k \big(\gb \tilde \omega_{x_i} + h' \big) \,.
\end{equation}
Introducing the white noise $W(\cdot)$ on $\R^d$ and considering the parallelepiped 
$\Delta := (-\frac{\delta^{a_1}}{2}, \frac{\delta^{a_1}}{2})\times\cdots \times  (-\frac{\delta^{a_d} }{2}, \frac{\delta^{a_{d}}}{2}) $,
we can replace each $\tilde\omega_x$ by
$\delta^{-d_{\eff}/2}W(x + \Delta)$. Since $h' = h' \delta^{-d_{\eff}} Leb(x+\Delta)$,
the inner sum in \eqref{eq:Zstart3} coincides
(recalling that $\psi^{(k)}_{\Omega_\delta}(x_1, \ldots, x_k)$ is piecewise constant),
with the following (deterministic + stochastic) integral:
\begin{equation} \label{eq:Zstart4}
	\idotsint_{\Omega^k}
	\psi^{(k)}_{\Omega_\delta}(x_1, \ldots, x_k)
	\, \prod_{i=1}^k \big( \gb \, \delta^{-d_{\eff}/2} \, W(\dd x_i)
        \,+\, h' \, \delta^{-d_{\eff}} \, \dd x_i  \big) \,.
\end{equation}
\smallskip

Finally, applying \eqref{eq:polyconv} and rescaling $\gb, h$ as in 
\eqref{eq:scalingpar},
equations \eqref{eq:Zstart3}-\eqref{eq:Zstart4} suggest that as $\delta \downarrow 0$
$Z^\omega_{\Omega_\delta; \gb, h}$ converges in distribution
to the continuum partition function $\bZ^W_{\Omega; \hbeta, \hh}$ defined
in \eqref{eq:Zcont}. The rigorous justification of these steps follows by
Theorems~2.3 and~2.5 in \cite{CSZ13}.
\end{proof}

\subsection{Some examples.} Let us look at some examples of application of Theorem \ref{th:mainhere}.
\vskip 2mm
{\bf The disordered pinning model.} \
Let $\tau = (\tau_k)_{k\ge 0}$ be a \emph{renewal process} on $\N$
with $\p(\tau_1 = n) = L(n) n^{-(1+\alpha)} $, $\alpha >0$ and $L(\cdot)$ a slowly varying function.
Consider $\Omega = (0,1)$,
$\delta = N^{-1}$ for $N \in \N$
and define $\P_{\Omega_\delta}^\re$ to be the law of
$(\sigma_x := \ind_{\delta\tau}(x))_{x\in \Omega_\delta}$,
where $\delta \tau = \{N^{-1} \, \tau_n\}_{n\ge 0}$ is viewed as a random subset of $\Omega$.
It is known that this model goes through a localization-delocalization transition when $h$ crosses a critical
value $h_c(\gb)$, \cite{G07}. 
If $\mathrm{M}(\gb)$ is the log-moment generating function of $\go$, then the quantity $h_c(\gb)+\mathrm{M}(\gb)$ equals zero
for sufficiently small $\gb$, if $\ga\in [0,1/2)$ (disorder irrelevance, \cite{A08}), 
is strictly positive for all $\gb>0$, if $\ga>1/2$, (disorder relevance, \cite{AZ09, DGLT09} ), while 
either scenarios are possible when $\ga=1/2$ (marginal), according to whether the quantity $\sum_{n=1}^N1/nL(n)^2$ diverges or not \cite{BL15}.
For $N\in\N$, $\hat\beta > 0$ and $\hat h\in\R$, set $h':=h+\mathrm{M}(\gb)$ and 
\begin{equation} \label{eq:scalingbetah}
   \begin{split}
    \beta_N = \begin{cases}
    \displaystyle
    \hbeta \frac{L(N)}{N^{\alpha - 1/2}} & \text{if } \ \frac{1}{2} < \alpha < 1 \\
    \displaystyle\rule{0em}{1.8em}\hbeta \frac{1}{\sqrt{N}} & \text{if } \ \E[\tau_1]<\infty
    \end{cases} \,,
    \qquad h'_N=
    \begin{cases}
    \displaystyle \hh \frac{L(N)}{N^{\alpha}} & \text{if } \ \frac{1}{2} < \alpha < 1 \\
	\displaystyle\rule{0em}{1.8em}\displaystyle \hh \frac{1}{N} & \text{if } \
	\E[\tau_1]<\infty
    \end{cases} \,.
  \end{split}
\end{equation}
As a consequence of Theorem \ref{th:mainhere}, we have that, when $\ga>1/2$,  the partition function
$Z_{Nt,\beta_N,h_N}^{\omega}$
of the disordered pinning model converges in distribution, for every $t \ge 0$, when $N\to\infty$
to the random variable $\bZ_{t,\hat\beta,\hat h}^{W}$ given by
\begin{equation} \label{eq:Zlimpin}
\bZ_{t,\hat\beta,\hat h}^{W} :=
1+\sum_{k=1}^\infty \frac{1}{k!}
\idotsint_{[0,t]^k}  \bpsi(t_1,\ldots,t_k) \,
\prod_{i=1}^k \big( \hbeta\,W(\dd t_i) + \hh\, \dd t_i \big) \,,
\end{equation}
where $W(\cdot)$ denotes white noise on $\R$
and $\bpsi_t (t_1,\ldots,t_k)$ is a symmetric function, defined for
$0<t_1<\cdots<t_k < t$ by
\begin{equation} \label{eq:cokec}
\bpsi(t_1,\ldots,t_k) =
\begin{cases}
\displaystyle\frac{C_\alpha^k }{t_1^{1-\alpha} (t_2-t_1)^{1-\alpha}
\cdots (t_k-t_{k-1})^{1-\alpha}} \qquad & \text{if } \ \frac{1}{2} < \alpha < 1 \\
\displaystyle\rule{0pt}{1.8em}
\frac{1}{\e[\tau_1]^k}  \qquad\qquad\qquad\qquad & \mbox{if } \ \E[\tau_1]<\infty
\end{cases} \,,
\end{equation}
with $C_\alpha := \frac{\alpha \sin(\pi\alpha)}{\pi}$.
The series in \eqref{eq:Zlimpin} converges in $L^2$, and in addition
$\bbE[(Z_{Nt,\beta_N,h_N}^{\omega})^2]
\to \bbE[(\bZ_{t,\hat\beta,\hat h}^{W})^2]$ as $N\to\infty$.
When $\E[\tau_1]<\infty$, e.g. when $\ga>1$, the continuum partition function is given concretely by 
\begin{equation}\label{pinWCag1}
 \bZ_{t,\hat\beta,\hat h}^{W} \overset{(d)}{=}
\exp\bigg\{\frac{\hbeta}{\e[\tau_1]}W_t+ \bigg(\frac{\hh}{\e[\tau_1]} - \frac{\hbeta^2}{2\e[\tau_1]^2}\bigg)t\bigg\} \,,
\end{equation}
where $W = (W_t)_{t\ge 0}$ denotes a standard Brownian motion.

\vskip 2mm
{\bf The directed polymer model}. \
Let $(S_n)_{n\ge 0}$ be a random walk on $\Z^d$. 
Consider $\Omega = (0,1) \times \R^d$, $\delta = N^{-1}$ for $N \in \N$
and set $\Omega_\delta := \big((\delta \Z) \times (\delta^{1/2}\Z^d)\big) \cap \Omega$.
Define $\P_{\Omega_\delta}^\re$ to be the law of the field
$(\sigma_x := \ind_{A_\delta}(x))_{x\in \Omega_\delta}$,
where $A_\delta := \{(\frac{n}{N}, \frac{S_n}{N^{1/\alpha}})\}_{n \ge 0}$ is viewed
as a random subset of $\Omega$. In this setting  the external field $h$ in \eqref{eq:model} 
 is taken to be equal $-\mathrm{M}(\gb)$.

For $d\geq3$ it is known that a small amount of disorder does not
change the $\sqrt{N}$ fluctuations of the simple, symmetric random walk, while superdiffusivity is expected and
localization occurs in $d=1,2$ for any arbitrary $\beta>0$, see \cite{CSY04} and references therein. When $d=1$, 
 $N^{2/3}$ fluctuations for the polymer path are predicted by the KPZ theory. In $d=2$ the fluctuation exponent is still elusive.
In $d=1$ we introduced  \cite{CSZ13}  a generalization of the directed polymer model where the increments of the walk lie
in the domain of attraction of an $\alpha$-stable law, with $1 < \alpha \le 2$.
 For $N\in\N$ and $\hat\beta > 0$ we set
$\beta_N := N^{-\frac{\alpha-1}{2\alpha}} \hat\beta $ and
the application of Theorem \ref{th:mainhere} in this situation yields that,
for every $t \ge 0$ and $x\in\R$, the point-to-point partition function
$Z_{Nt,\beta_N}^{\omega}(N^{1/\alpha}x)$ converges in distribution as $N\to\infty$ to
a random variable $\bZ_{t,\hat\beta}^{W}(x)$, which is a solution to the fractional Stochastic Heat Equation
\begin{equation}\label{eq:SFHE}
\partial_t  u =\Delta^{\alpha/2} \, u +\sqrt{p}\hat\beta \, \dot W\,u \qquad \text{with}\qquad 
u(0, \cdot ) =\delta_0(\cdot) ,
\end{equation}
where $p$ is the period of the walk.
When $\ga=2$ the fractional Laplacian is the usual Laplacian on $\bbR$ and this result was obtained in \cite{AKQ14}.
\vskip 2mm
{\bf The random field Ising model $(d=2)$.} \
Consider a bounded and connected set $\Omega \subseteq \R^2$ with smooth boundary
and define $\P_{\Omega_\delta}^\re$ to be the \emph{critical} Ising model on $\Omega_\delta$
at inverse temperature $\beta = \beta_c =
\frac{1}{2}\log(1+\sqrt{2})$ and $+$ boundary condition.
The measure on spin configurations $\{\pm 1\}^{\Omega_\gd}$ is given by
\begin{align*}
\P^{\re}_{\Omega_\gd}(\gs):=\P^+_{\Omega_\gd}(\gs)=\frac{1}{Z^+_{\Omega_\gd}} \exp\Big(\gb_c\sum_{x\sim y\in \Omega_\gd \cup \partial \Omega_\gd} \gs_x\gs_y\Big)
\prod_{x\in\partial \Omega_\gd} \ind_{\gs_x=+1}
\end{align*}
It has been shown in \cite{CHI12} that Assumption \ref{ass:main} is satisfied for $\P^+_{\Omega_\gd}$ with $\gamma=1/8$. Theorem \ref{th:mainhere}
allows to construct the continuum partition function of the random field perturbation of $\P^+_{\Omega_\gd}$. The continuum partition function is given in terms of iterated Wiener integrals but it is less explicit \cite{CSZ13}.
 \begin{remark}{\rm
 Theorem \ref{th:mainhere} allows to construct the continuum limit of partition functions. Even though this is a main ingredient, 
 additional work is required in order to construct the
 continuum limit of the disordered measure. This is related to establishing convergence of partition functions as a process.
 In this way, the construction of the continuum models has been achieved for the pinning model \cite{CSZ16} and for the directed polymer model \cite{AKQ14b} but not yet for the random field Ising model. 
 }
 \end{remark}
\section{Marginal relevance}
Inspecting the first part of relation \eqref{eq:scalingpar}, we see that the situation $\gamma=d_\eff/2$ is incompatible with the assumption of 
{\it weak disorder} as the
exponent in the scaling vanishes and the strength $\gb$ of the disorder does not tend to zero when $\delta\downarrow 0$. 
This is essentially a restatement of Harris' marginal condition $\nu=2/d$
(see \cite{CSZ13} for a discussion). In the case of the pinning model the incarnation of this condition is $\alpha=1/2$ 
and in the case of the long range DPRM $\alpha=1$, which were
both excluded. The issue with these marginal values is not
 technical but rather structural. This can been seen if one sets, for example, the marginal value
$\alpha=1/2$ in the formula for the continuum pinning partition \eqref{eq:cokec}. This formal substitution leads to stochastic integrals, which cannot be given an It\^o sense
due to the fact that the kernel \eqref{eq:cokec} is not $L^2$ integrable. The situation is similar in the case of $(2+1)$-dimensional
 directed polymer corresponding to a simple symmetric random walk, since
 $d_\eff=2$, while by the local limit theorem, the exponent $\gamma$ in \eqref{eq:polyconv} takes the value $d_{\eff}/2=1$.
 
 A source of this difficulty is that the vanishing of the exponent $d/2-\gamma$ induces a new time scale, which is actually exponential. To 
 be more
 precise, let us look at the multilinear expansion \eqref{eq:Zstart2}. In this setting we will be assuming that $\bbE[\eta_x]=0$ and,
for the sake of exposition,  that the $k-$point correlation function 
 is the correlation function of either a simple, symmetric 
 random walk on $\bbZ^2$, or a one dimensional Cauchy random walk, i.e. 
 $\psi((n_1,x_1),...,(n_k,x_k))=\P(S_{n_1}=x_1,...,S_{n_k}=x_k)=\prod_{i=1}^n\P(S_{n_i-n_{i-1}}=x_i-x_{i-1})$. 
 We will denote by $q_n(x)$ the probability that the random walk is at location $x$ at time $n$.
  Computing the variance of the $k$-th term in the multilinear expansion,
  we obtain that it is asymptotically equal to 
  \begin{align*}
  \bbE\left[\Big(\sumtwo{x_1,...,x_k}{1\leq n_1<\cdots<n_k\leq N} \,\prod_{i=1}^k q_{n_i-n_{i-1}}(x_i-x_{i-1}) \prod_{i=1}^N \eta_{x_i} \Big)^2\right] 
  &=\sumtwo{x_1,...,x_k}{1\leq n_1<\cdots< n_k\leq N} \, \prod_{i=1}^k q_{n_i-n_{i-1}}(x_i-x_{i-1})^2 \\
  &\approx (\log N)^k.
  \end{align*}
  This computation indicates two things: First, that the right rescaling of $\gb$ is $\gb_N:=\hbeta/\sqrt{\log N}$. Second,
  that the asymptotic variance will remain unchanged if we summed over a time horizon $tN$ for any arbitrary but fixed time variable $t$. On the other hand,
  if we considered a time horizon $N^t$ with $t>0$, we will see a 
  (linear in $t$) change in the asymptotic behavior. 
 Therefore, the correct time scale is $(N^t\colon t>0)$ and to incorporate this we decompose the summations over $n_1,...,n_k$ over intervals
 $n_j-n_{j-1}\in I_{i_j}$, with $I_{i_j}=\big(N^{\frac{i_j-1}{M}}, N^{\frac{i_j}{M}} \big]$, $i_j\in\{1,...,M\}$ and $M$ a coarse graining parameter which will eventually
 tend to infinity. We can then rewrite the $k$-th term in the expansion \eqref{eq:Zstart2} as  
\begin{align}\label{thetas}
&\frac{\hat\beta^k}{M^{k/2}}\sum_{1\leq i_1,...,i_k\leq M}  \Theta^{N,M}_{i_1,...,i_k} \qquad \text{where}\\
&\Theta^{N,M}_{i_1,...,i_k} :=
\left(\frac{M}{\log N}\right)^{k/2}\sumtwo{n_j-n_{j-1}\in\, I_{i_j}\,\,\text{for }\,j=1,...,k}{x_1,...,x_k\in \mathbb{Z}^d}\,\, \prod_{j=1}^k q_{n_j-n_{j-1}}(x_j-x_{j-1}) \,\eta_{(n_j,z_j)}
\notag
\end{align}
We now observe that if an index $i_j$ is a running maximum for the $k$-tuple $\bi:=(i_1,...,i_k)$, i.e. $i_j> \max\{i_1,...,i_{j-1}\}$ 
then $\big(N^{\frac{i_j-1}{M}}, N^{\frac{i_j}{M}} \big] \ni n_j \gg 
 n_{r} \in \big(N^{\frac{i_r-1}{M}}, N^{\frac{i_{r}}{M} }\big] $, for all
$r<j$
\footnote{Strictly speaking, for this inequality to be valid uniformly, we need to restrict to values of
 $\bi\in \{1,...,M\}^k_\sharp:=\{\bi \in \{1,...,M\}^k \colon |i_j-i_k|>1\}$,
 but this is a minor technical point that be easily taken care of.}.
Using a version of the local limit theorem, the inequality implies that
$ q_{n_j-n_{j-1}}(z_j-z_{j-1}) \approx q_{n_{j}}(z_{j}) $ for $n_j\in I_{i_j}$ and $n_{j-1}\in I_{i_{j-1}}$. 
Decomposing the sequence $\bi:=(i_1,...,i_k)$  according to its running maxima, i.e.
 $\bi=(\bi^{(1)},...,\bi^{(\mathfrak{m}(\bi))})$ with $\bi^{(r)}:=(i_{\ell_{r}},...,i_{\ell_{r+1}}-1)$ and 
 $i_{1}=i_{\ell_1}<i_{\ell_2}<\cdots<i_{\ell_{\mathfrak{m}}}$ the successive running maxima, it can be shown 
 that  \eqref{thetas} factorizes (asymptotically when $N$ tends to infinity) as
\begin{align*}
 \frac{\hbeta^k}{M^{\frac{k}{2}}} \sum_{\bi\in \{1,\ldots, M\}^k_\sharp}
\Theta^{N;M}_{\bi^{(1)}} \Theta^{N;M}_{\bi^{(2)}}\cdots \Theta^{N;M}_{\bi^{(\mathfrak{m})}}.
\end{align*}
The heart of the argument is to show that all the $\Theta^{N;M}_{\bi^{(j)}}$ converge jointly,
 when $N\to\infty$ to standard normal variables. This is established
in Proposition 5.2 in \cite{CSZ15}. The argument is combinatorial and makes use of (a version of) the Fourth Moment Theorem \cite{dJ90, NP05, NPR10}. Once this 
convergence is established, a re-summation leads to   
\begin{theorem}[Limit of partition functions \cite{CSZ15}]\label{thm:subcritical}
Let $Z_{N,\beta}^\omega$ be the partition function of
a directed polymer or a pinning model whose transition probability kernel or renewal function satisfies a form of local limit theorem.
Assume that the replica overlap of two independent random walks $S,S'$  
(in case of directed polymers) or two independent renewals $\tau,\tau'$ (in case of pinning)
\begin{align*}
\Ro_N :=\begin{cases}
\E\big[\,\sum_{n=1}^N \ind_{S_n=S'_n}\,\big],\\
\\
 \E\big[\,\sum_{n=1}^N \ind_{n\in\tau\cap \tau'}\,\big],
 \end{cases}
\end{align*}
diverges as a slowly varying function when $N\to\infty$. Then, defining
$\beta_N := \hat\beta/\sqrt{\Ro_N} $ with $ \hat\beta \in (0,\infty)$,
the following convergence in distribution holds,
where $(W_t)_{t\ge 0}$ is a standard Brownian motion:
\begin{equation} \label{eq:conve0}
	Z_{N,\beta_N}^\omega \xrightarrow[\,N\to\infty\,]{(d)} 
	\bs{Z}_{\hat\beta} :=
	\begin{cases}
	\exp\bigg( \displaystyle\int_0^1\textstyle
	\frac{\hat\gb }{\sqrt{1-\hat\gb^2 \,t}} \, \dd W_t-
	\frac{1}{2} \displaystyle\int_0^1\textstyle
	\frac{\hat\gb^2\ }{1-\hat\gb^2\,t} \, \dd t \bigg)
	& \text{if } \hat \beta < 1 \\
	0 & \text{if } \hat\beta \ge 1
	\end{cases} \,.
\end{equation}
Moreover, for $\hat\beta < 1$ one has
$\lim_{N\to\infty}\bbE[(Z_{N,\beta_N}^\omega)^2] 
= \bbE[(\bs{Z}_{\hat\beta})^2]$.
\end{theorem}
An extension of the above Theorem at a process level
 is also established in Theorems 2.11 and 2.12 in \cite{CSZ15} 
 which show a multi scale and log-correlated structure. 
 The analogous results for the $2d$ Stochstic Heat Equation are also established in \cite{CSZ15}. 

It is worth pointing out the phase transition at $\hbeta=1$. It is by now well known \cite{CSY04} that there are two regimes in 
 directed polymers in random media: the strong and the weak disorder regimes, which can be classified as the regimes where
 the normalized partition function $ Z_{N,\beta}^\go / \bbE[Z_{N,\beta}^\go]$ has an a.s. zero or strictly positive limit, as $N$ tends to infinity.
 In the weak disorder regime the polymer path behaves, essentially, as a simple symmetric 
 random walk while in the strong disorder regime localization takes place.
 It was shown \cite{CSY04} that for $d\geq 3$ and $\beta$ smaller than a critical value $\beta_c$, the model is in the weak disorder regime.
 On the other hand in dimensions $d=1,2$ it is shown \cite{CSY04} that the critical value $\beta_c$ equals zero, that is, for any $\beta>0$
 there is strong disorder. Theorem \ref{thm:subcritical} shows that actually there is a transition 
 between weak and strong disorder in dimension two, which is observed when one looks on the scale $\hat\beta/\Ro_N$ and
 with transition at $\hbeta_c=1$. Moreover, this transition is common among all marginally relevant disordered polymer models that fall 
 within the scope of Theorem \ref{thm:subcritical}. Finally, it is worth noting that similar transition takes place in the model of Gaussian 
 Multiplicative Chaos \cite{RV13}.

\section{Some consequences}

Having a continuum model at hand can be very useful in extracting sharp information on the phase diagram of the discrete models. To describe one such use, 
 we will use the pinning 
model as an example. Its free energy (we recenter the parameter $h$ to $h-\mathrm{M}(\gb)$ for convenience) is defined as 
\begin{align*}
f(\beta, h):=\lim_{N\to \infty} \frac{1}{N}\log Z^\go_{N,\beta,h-\mathrm{M}(\gb)} = \lim_{N\to \infty} \frac{1}{N} \bbE \big[\log Z^\go_{N,\beta,h-\mathrm{M}(\gb)}\big].
\end{align*}
One can show that $f$ is nonnegative and transition from a zero to a strictly positive value determines a localization/delocalization transition. In particular,
we can define the critical value $h_c(\beta):=\sup\{ h\colon F(\beta,h)=0 \}$. Much effort has been devoted to the study of this critical value and in particular whether
it is strictly positive for all sufficiently small values of $\beta$. It has been 
shown that $h_c(\beta)=0$, for all sufficiently small $\beta$, when $\ga\in [0,1/2)$ \cite{A08, T08}, 
while if $\ga\in(1/2,1)$ there exist constants $c,C$ and a slowly varying function $\tilde{L}(\cdot)$ determined by $L(\cdot)$, such that 
$c\beta^{\frac{2\alpha}{2\alpha-1}}\tilde{L}(1/\beta^2) < h_c(\beta) < C\beta^{\frac{2\alpha}{2\alpha-1}}\tilde{L}(1/\beta^2)$, for all sufficiently
small $\gb$ \cite{AZ09, DGLT09}. With the help of a continuum partition function, 
 this result can be strengthened to a statement of existence of the limit  
 $ \lim_{\gb\to0}\gb^{-\frac{2\ga}{2\ga-1}}\tilde{L}(1/\beta^2)^{-1}  h_c(\gb)$. 
Denote the continuum free energy (which we assume to exist) by 
\begin{align*}
{\bf F}(\hbeta,\hat h)=\lim_{t\to\infty} \frac{1}{t} \bbE [\log \bZ_{t,\hbeta, \hat h}] =\lim_{t\to\infty} \frac{1}{t} \lim_{N\to\infty} \bbE \log Z_{Nt,\gb_N,h_N-\mathrm{M}(\gb_N)},
\end{align*}
and assuming that we can interchange the limits over $N,t$, this leads to
\begin{align*}
{\bf F}(\hbeta,\hat h)=\lim_{N\to\infty} N \lim_{t\to\infty} \frac{1}{Nt} \bbE \log Z_{Nt,\gb_N,h_N-\mathrm{M}(\gb_N)} = \lim_{N\to\infty} N F(\gb_N,h_N). 
\end{align*}
Expressing $N$ and $h$ in terms of $\gb$, we can write the above as
\begin{align}\label{limits}
{\bf F}(\hbeta,\hat h)=\lim_{\beta\to0} \frac{F(\gb, \gb^{\frac{2\ga}{2\ga-1}} \tilde{L}(\frac{1}{\gb}) \hat h)}{ \gb^{\frac{2}{2\ga-1}} \,\tilde{L}(\frac{1}{\gb})^2},
\,\,\,\text{indicating that}\,\,\,
{\bf h}_c(\hbeta) =\lim_{\gb\to 0} \frac{h_c(\gb)}{\gb^{\frac{2\ga}{2\ga-1}} \tilde{L}(\frac{1}{\gb})}.
\end{align}
The quantitative estimate that justifies the above interchange 
of limits and thus the validity of the limits in \eqref{limits} was shown in \cite{CTT15} and is
\begin{theorem}[\cite{CTT15}]
Let $F(\gb,h)$ be the free energy of the disordered pinning model with renewal exponent $\ga\in(1/2,1)$ and 
${\bf F}(\hbeta,\hat h)$ the free energy of the corresponding continuum model. For all $\hbeta>0, \hat h\in \bbR$ and 
$\eta>0$ there exist $\gd_0>0$ such that for all $\gd\in (0,\gd_0)$
\begin{align*}
{\bf F}(\hbeta,\hat h-\eta)\leq \frac{F(\hbeta\, \gd^{\ga-\frac{1}{2}} L(1/\gd),  \hat h\, \gd^\ga L(1/\gd))}{\gd}
\leq {\bf F}(\hbeta,\hat h+\eta).
\end{align*}
\end{theorem}
The situation should be similar in other cases, like the Directed Polymer and the Random Field Ising model but the corresponding quantitative estimates have not been established, yet. We expect that one should be able to prove that for the $1d$ DPRM the limit 
$\lim_{\gb\to0}\gb^{-4} f(\gb)$ should exist
and the same should be the situation for the limit $\lim_{h\to0}h^{-1/15} \langle\sigma_0\rangle_{\gb_c,h}$ for the $2d$ RFIM at critical 
temperature, with the limits being given in terms of the 
corresponding continuum free energies.
\begin{center}{\bf Acknowledgements}
\end{center}
 FC acknowledges the support of GNAMPA-INdAM. RS is supported by NUS grant R-146-000-185-112. NZ is funded by EPSRC via EP/L012154/1.


\begin{thebibliography}{10}

\bibitem{AW90}
M.~Aizenman and J.~Wehr,
Rounding effects of quenched randomness on first-order phase transitions,
{\em Commun.\ Math.\ Phys.} 130, 489--528, 1990

\bibitem{AKQ14}
T. ~Alberts, K.~ Khanin, J.~ Quastel,
Intermediate Disorder for $1+1$ Dimensional Directed Polymers,
{\em Ann. Probab.} 42, 1212--1256, 2014.

\bibitem{AKQ14b}
T. Alberts, K. Khanin, J.Quastel,
The continuum directed random polymer,
 J. Stat. Phys. 154.1-2 (2014) 305-326

\bibitem{A08}
K.~S.~Alexander,
The effect of disorder on polymer depinning transitions,
{\em Commun. Math. Phys. 279 (2008), no. 1, 117–146}

\bibitem{AZ09}
K.~S.~Alexander and N.~Zygouras.
Quenched and annealed critical points in polymer pinning models.
\textit{Commun. Math. Phys.} 291 (2009), no. 3, 659–689.

\bibitem{AKQ14a}
T.~Alberts, K.~Khanin, and J.~Quastel.
The intermediate disorder regime for directed polymers in dimension $1+1$.
{\em Ann.\ Probab.} 42 (2014), 1212--1256.

\bibitem{BL15}
 Q. Berger, H. Lacoin, 
 Pinning on a defect line: characterization of marginal disorder relevance and sharp asymptotics for the critical point shift, 
{\em Journal of the Inst. Math. Jussieu}, (2015) to appear

\bibitem{BCPSZ14}
Q.~Berger, F.~Caravenna, J.~Poisat, R.~Sun, and N.~Zygouras.
The Critical Curve of the Random Pinning and Copolymer Models at Weak Coupling.
{\em Commun. Math. Phys.} 326 (2014),  507--530.

\bibitem{BK88}
J.~Bricmont and A.~Kupiainen.
Phase transition in the 3d random field Ising model.
{\em Commun.\ Math.\ Phys.} 116 (1988), 539--572.

\bibitem{CGN12b}
F.~Camia, C.~Garban, and C.M.~Newman.
The Ising magnetization exponent on $\Z^2$ is $1/15$.
{\em Probab. Theory Relat. Fields} (to appear), arXiv:1205.6612.

\bibitem{CSZ13}
F. Caravenna, R. Sun, N. Zygouras.
Polynomial chaos and scaling limits of disordered systems.
{\em J. Eur. Math. Soc.}, to appear.

\bibitem{CSZ15}
F. Caravenna, R. Sun, N. Zygouras.
Universality in marginally relevant disordered systems,
arXiv:1510.06287 

\bibitem{CSZ16}
F.~Caravenna, R.~Sun, N.~Zygouras.
The continuum disordered pinning model.
{\em Prob. Theory Rel. Fields}, (2016), no. 1-2, 17-59

\bibitem{CTT15}
F. Caravenna, F.L. Toninelli, N. Torri, 
Universality for the pinning model in the weak coupling regime,
 arXiv:1505.04927

\bibitem{CCFS86}
J.T. Chayes, L. Chayes, D.S. Fischer, and T. Spencer.
Finite-Size Scaling and Correlation Lengths for Disordered Systems.
\emph{Phys. Rev. Lett.} 57, number 24, 15 December 1986.

\bibitem{CCFS89}
J.T. Chayes, L. Chayes, D.S. Fischer, and T. Spencer.
Correlation Length Bounds for Disordered Ising Ferromagnets.
\emph{Commun. Math. Phys.} 120 (1989), 501--523.

\bibitem{CHI12}
D.~Chelkak, C.~Hongler, and C.~Izyurov.
Conformal invariance of spin correlations in the planar Ising model.
{\em Ann. Math.} (to appear), arXiv:1202.2838.

\bibitem{CSY04}
F. Comets, T. Shiga, N. Yoshida.
Probabilistic analysis of directed polymers in a random environment: a review.
{\em Stochastic analysis on large scale interacting systems}, 115--142, 
Adv. Stud. Pure Math., 39, Math. Soc. Japan, Tokyo, 2004.

\bibitem{dJ90}
P. de Jong,
A central limit theorem for generalized multilinear forms,
{\em J. Multivariate Anal.} 34,  275--289, 1990

\bibitem{DGLT09}
B. Derrida, G. Giacomin, H. Lacoin, F. Toninelli, 
Fractional moment bounds and disorder relevance for pinning models.
{\em Comm. Math. Phys.} 287.3 (2009),  867-887

\bibitem{G07}
G.~Giacomin.
Random Polymer Models.
{\em Imperial College Press}, London, 2007.

\bibitem{G10}
G.~ Giacomin.
Disorder and critical phenomena through basic probability models.
\'Ecole d'\'Et\'e de Probabilit\'es de Saint-Flour XL – 2010.
{\em Springer Lecture Notes in Mathematics} 2025.

\bibitem{GT06}
G.~ Giacomin, F. Toninelli,
Smoothing effect of quenched disorder on polymer depinning transitions,
{\em Comm. Math. Phys.} 266.1 (2006) 1-16

\bibitem{H74}
A.B.~Harris,
Effect of random defects on the critical behaviour of Ising models.
{\em Journal of Physics C: Solid State Physics}, 1974.


\bibitem{L10}
H.~Lacoin.
New bounds for the free energy of directed polymer in dimension 1+1 and 1+2.
{\em Commun. Math. Phys.} 294 (2010) 471--503.


\bibitem{MOO10}
E.~Mossel, R.~O'Donnell, and K.~Oleszkiewicz.
Noise stability of functions with low influences: Variance and optimality.
{\em Ann.\ Math.} 171, 295--341, 2010.

\bibitem{NPR10}
I. Nourdin, G. Peccati, G. Reinert.
Invariance principles for homogeneous sums: universality of Gaussian Wiener chaos.
{\it Ann. Probab.} 38, 1947--1985, 2010.

\bibitem{NP05}
D. Nualart and G. Pecati.
Central limit theorems for sequences of multiple stochastic integrals.
{\it Ann. Probab.} 33, 177--193, 2005.

\bibitem{RV13}
R. Rhodes and V. Vargas,
Gaussian multiplicative chaos and applications: a review.
{\em Probab. Surveys} 11, 315--392, 2014.

\bibitem{T08}
F.L. Toninelli
A replica-coupling approach to disordered pinning models, 
{\em Comm. Math. Phys.} 280, 389-401 (2008)
\end{thebibliography}
\end{document}